\documentclass[twoside]{amsart}

\usepackage{latexsym,amssymb}

\def\demo{\noindent{\bf Proof. }}
\def\sqr#1#2{{\vcenter{\hrule height.#2pt
        \hbox{\vrule width.#2pt height#1pt \kern#1pt
                \vrule width.#2pt}
        \hrule height.#2pt}}}
\def\square{\mathchoice\sqr64\sqr64\sqr{4}3\sqr{3}3}
\def\QED{\hfill$\square$}

\def\tratto{\mbox{\rule{2mm}{.2mm}$\;\!$}}

\def\frak{\mathfrak}

\newtheorem{Theorem}{Theorem}[section]

\newtheorem{Corollary}[Theorem]{Corollary}

\newtheorem{Remark}[Theorem]{Remark}

\newtheorem{Definition}[Theorem]{Definition}

\begin{document}

\baselineskip=16pt

\title[Residually $S_2$ ideals and projective dimension one modules]
{{\Large\bf On Residually $S_2$ ideals and \\
projective dimension one modules}}

\author[A. Corso and C. Polini]
{Alberto Corso \and Claudia Polini}

\thanks{
AMS 1991 {\it Mathematics Subject Classification}. Primary 13A30;
Secondary 13B22,  13C10, 13C40. \newline\indent
\indent
{\bf Acknowledgements:} Both authors sincerely thank Bernd Ulrich
for many helpful discussions they had concerning the material in this paper.
The NSF, under grant DMS-9970344, has also partially supported the
research of the second author and has therefore her heartfelt thanks.}

\address{Department of Mathematics, Michigan State University, E.
Lansing, Michigan 48824}
\email{corso@math.msu.edu}

\address{Department of Mathematics, Hope College, Holland, Michigan 49422}
\email{polini@cs.hope.edu}

\begin{abstract}
We prove that certain modules are faithful. This enables us to draw consequences about
the reduction number and the integral closure of some classes of ideals.
\end{abstract}

\maketitle

\section{Introduction}

This work has been motivated by a recent paper of
C. Huneke about cancellation theorems for special ideals in local Gorenstein
rings. In its simplest form any such theorem says that if $I, J$ and $L$ are
ideals in a Noetherian local ring $R$ such that $LI \subset JI$ then
$L \subset J$. Using the so called determinant trick, in general one can
only conclude that $L \subset JI \colon I \subset \overline{J}$,
where $\overline{J}$ denotes the integral closure of $J$.
We recall that the {\it integral closure} $\overline{J}$ of an
ideal $J$ is the set $($ideal, to be precise$)$ of all elements integral over
$J$. An element $x \in R$ is {\it integral} over $J$ if $x$ satisfies a monic
equation of the form $x^n+a_1 x^{n-1} + \ldots + a_{n-1} x + a_n =0$,
where $a_i \in J^i$ for $ 1 \leq i \leq n$.

To say that cancellation holds for every such ideal $L$ is equivalent to say
that the $R/J$-module $I/JI$ is faithful, that is to say $JI \colon I =
J$.
There aren't many known instances when cancellation happens. Besides
the ones described in \cite{H}, we would like to single out
\cite[Theorem~2.4 and Theorem~3.1]{CHV}.

More generally, one can ask for which ideals $I$ and $J$ and which positive
integers $t$ the equality $JI^t \colon I^t = J$ holds.
This kind of question is particularly interesting when the ideal $J$ is a
$($minimal$)$ reduction of the ideal $I$. We say that a subideal $J$ of $I$
is a {\it reduction} of $I$ if $I^{r+1} = JI^r$ for some positive integer
$r$. The smallest such $r$ is called the {\it reduction number} of $I$ with respect
to $J$. {\it Minimal reductions} are reductions minimal with respect to containment.
If the residue field is infinite, their minimal number of generators does
not depend on the minimal reduction of the ideal. This number is called
{\it analytic spread} of $I$, in symbols $\ell=\ell(I)$. It equals the dimension
of the fiber cone of $I$ and it is always greater than or equal to the
height of $I$.
If $I^{r+1}=JI^r$ one has that $I \subset JI^r \colon I^r$. Thus, conditions
of the kind $JI^t \colon I^t = J$, with $J$ a reduction of $I$ and $t$ some
integer, give a severe bound on the reduction number of ideals which admit
proper reductions.

Our main results, in the ideal case, are Theorem~\ref{thm1} and
Theorem~\ref{thm2}. They answer the above questions in the case of ideals
satisfying local bounds on the minimal number of generators up to a given
codimension and some other residual properties in the sense of
\cite{CEU, CPU, JU, PU, PU2, U, U2}.
In the same spirit as \cite{H}, immediate consequences of Theorem~\ref{thm2}
are drawn about the reduction number and the integral closure of an ideal $I$.
The proofs of Theorem~\ref{thm1} and Theorem~\ref{thm2} use, though,
techniques coming essentially from \cite{U} $($and its subsequent
refinements$)$.

We also prove a similar cancellation theorem in the case of certain modules
of projective dimension one $($see Theorem~\ref{pd1}$)$.
The motivations and the techniques we use come from \cite{RAM}, where the
authors study Rees algebras of modules via Bourbaki ideals.

It is worth pointing out the analogy between the classes of ideals and
modules presently studied and the ones recently studied in \cite{CPU}:
some of the theorems proved in \cite{CPU} are also sort of cancellation
theorems.

\medskip

\section{The main results}

\subsection{The ideal case}

We first recall some additional definitions that are essential in stating
our results. For a more comprehensive treatment, we refer the reader to
\cite{U, CEU}.

\begin{Definition}
{\rm Let $R$ be a local Cohen--Macaulay ring, let $I$ be an $R$-ideal
of grade $g$, and let $s \geq g$ be an integer.

\begin{enumerate}
\item
We say that $I$ satisfies property $G_s$, if $\mu(I_{\frak p}) \leq
{\rm dim}\, R_{\frak p}$ for any prime ideal ${\frak p}$
containing $I$ with ${\rm dim}\, R_{\frak p} \leq s-1$.

\item
An $s$-residual intersection of $I$ is a proper ideal $K= {\mathfrak a}
\colon I$ where ${\mathfrak a}$ is a subideal of $I$ with $\mu({\mathfrak a})
\leq s \leq {\rm height}(K)$.

\item
An $s$-residual intersection $K$ of $I$ is called a geometric
$s$-residual intersection if ${\rm height}(I + K) \geq s+1$.

\item
We say that $I$ is $s$-residually $S_2$ if for every $g \leq i \leq s$ and
every $i$-residual intersection $K$ of $I$, $R/K$ satisfies Serre's
condition $S_2$.

\item
We say that $I$ is weakly $s$-residually $S_2$ if for every $g \leq i
\leq s$ and every geometric $i$-residual intersection $K$ of $I$,
$R/K$ satisfies Serre's condition $S_2$.
\end{enumerate} }
\end{Definition}

\begin{Theorem}\label{thm1}
Let $R$ be a local Cohen--Macaulay ring. Let $I$ be an unmixed $R$-ideal
of grade $g$ satisfying $G_s$ for some $s \geq g+1$.
Let $K = {\mathfrak a} \colon I$ be an $s$-residual intersection.
Then there exists a generating sequence $a_1, \ldots, a_s$ of ${\mathfrak a}$
such that with ${\mathfrak a}_i = (a_1, \ldots, a_i)$ then
\begin{itemize}
\item[$(${\it a}$)$]
${\mathfrak a}_i I^t \colon I^t = {\mathfrak a}_i$ if $1 \leq i \leq g$
and for every $t$;

\item[$(${\it b}$)$]
${\mathfrak a}_i I^t \colon I^t \subset I$ if $g \leq i \leq s$
and for every $t$.
\end{itemize}
Suppose further that $I$ is weakly $(s-1)$-residually $S_2$. Then
\begin{itemize}
\item[$(${\it c}$)$]
${\mathfrak a}_i I^t \colon I^t = {\mathfrak a}_i$ if $g+1 \leq i
\leq s-1$ or if $i=s$ and $K$ is a geometric $s$-residual intersection,
and for every $t$.
\end{itemize}
\end{Theorem}
\demo
Pick the elements $a_i$ as in \cite[Corollary~1.6$(${\it a}$)$]{U}.

$(${\it a}$)$
Clearly, ${\mathfrak a}_i \subset {\mathfrak a}_i I^t \colon I^t$ for all
$t$.
It is enough to prove the statement at the associated primes of $R/{\mathfrak
a}_i$ which all have height $i$. Let ${\mathfrak p}$ be any such prime. If
$i < g$ then $({\mathfrak a}_i)_{\mathfrak p} I^t_{\mathfrak p}
\colon I^t_{\mathfrak p} = ({\mathfrak a}_i)_{\mathfrak p} : R_{\mathfrak p} =
({\mathfrak a}_i)_{\mathfrak p}$. If $i=g$ we either have $I_{\mathfrak p} =
({\mathfrak a}_g)_{\mathfrak p}$ or $I_{\mathfrak p} = R_{\mathfrak p}$.
In both case we conclude that $({\mathfrak a}_g)_{\mathfrak p}
I^t_{\mathfrak p} \colon I^t_{\mathfrak p} = ({\mathfrak a}_g)_{\mathfrak p}$.

$(${\it b}$)$
It is enough to show the inclusion at the associated
primes of $I$, which all have height $g$ since $I$ is unmixed. Let
${\mathfrak p}$ be any such prime. Then $I_{\mathfrak p} =
({\mathfrak a}_g)_{\mathfrak p}$ and
$({\mathfrak a}_i)_{\mathfrak p} I^t_{\mathfrak p} \colon
I^t_{\mathfrak p} = ({\mathfrak a}_g)_{\mathfrak p}
= I_{\mathfrak p}$.

$(${\it c}$)$
Set $K_i = {\mathfrak a}_i \colon I$ and observe that
${\mathfrak a}_i I^t \colon I^t \subset {\mathfrak a}_i
\colon I^t = ({\mathfrak a}_i \colon I) \colon I^{t-1} = K_i \colon
I^{t-1}$. Since $K_i$ is a $i$-geometric
residual intersection and is unmixed of height $i$, it follows that
$I^{t-1}$ has positive grade modulo $K_i$. Thus, we conclude that
$K_i \colon I^{t-1} = K_i$. Using $(${\it b}$)$, we have
\[
{\mathfrak a}_i I^t \colon I^t \subset K_i \cap I = {\mathfrak a}_i,
\]
as desired. We observe that the equality $K_i \cap I = {\mathfrak a}_i$ follows
from \cite[Proposition~3.4$(${\it c}$)$]{CEU}. \QED

\medskip

Our next goal is to relax the assumptions of Theorem~\ref{thm1}$(${\it c}$)$
in the case $i=s$.

\begin{Theorem}\label{thm2}
Let $R$ be a local Gorenstein ring. Let $I$ be an $R$-ideal
of grade $g$, satisfying $G_s$ for some $s \geq g$ and locally unmixed
in codimension $s$. Suppose further that
$I$ is weakly $(s-1)$-residually $S_2$ and that
${\rm Ext}^{g+j}_{R}(R/I^j, R)_{\mathfrak p} = 0$
whenever ${\mathfrak p} \in V(I)$ with ${\rm dim}\, R_{\mathfrak p}\leq s$
and $1 \leq j \leq s-g$.
Let ${\mathfrak a} \colon I$ be an $s$-residual intersection.
Then for every $1 \leq t \leq s-g$
\[
{\mathfrak a} I^t \colon I^t = {\mathfrak a}.
\]
\end{Theorem}
\demo
Again, pick a generating sequence $a_1, \ldots, a_s$ of ${\mathfrak a}$
as in \cite[Corollary~1.6$(${\it a}$)$]{U} and we may assume $s \geq g+1$.
Clearly, ${\mathfrak a} \subset {\mathfrak a} I^t \colon I^t$ for all $t$.
Hence it is enough to show the equality at the associated primes of
$R/{\mathfrak a}$ which have height at most $s$ $($see
\cite[Proposition~3.4$(${\it b}$)$]{CEU}$)$. Let ${\mathfrak p}$ be any
such prime. We may also assume that ${\mathfrak p} \in V(I)$, as otherwise
the conclusion is immediate. Since ${\mathfrak a} \colon I$
is an $s$-residual intersection,
it follows that if ${\mathfrak p}$ has height at most $s-1$ then
$I_{\mathfrak p} = {\mathfrak a}_{\mathfrak p}$. The unmixedness assumption
on $I$ forces then ${\mathfrak p}$ to have height $g$.
Therefore $I_{\mathfrak p} =
{\mathfrak a}_{\mathfrak p} = ({\mathfrak a}_g)_{\mathfrak p}$, where
${\mathfrak a}_g = (a_1, \ldots, a_g)$, and
${\mathfrak a}_{\mathfrak p} I^t_{\mathfrak p} \colon
I^t_{\mathfrak p} = ({\mathfrak a}_g)_{\mathfrak p}
= {\mathfrak a}_{\mathfrak p}$.

In conclusion we may assume that ${\mathfrak p}$ has height $s$ and
therefore, by the unmixedness assumptions on $I$, we have that
$I_{\mathfrak p} \not= {\mathfrak a}_{\mathfrak p}$.
After a change of notation, let $R$ be a local Gorenstein ring of
dimension $s$ and $I$ an $R$-ideal of grade $g$ satisfying $G_s$ and
${\rm Ext}^{g+j}_R (R/I^j, R) = 0$ for $1 \leq j \leq s-g$.
Let `$\,{}^{\tratto}$' denote images modulo $K_{s-1} = {\mathfrak a}_{s-1}
\colon I$. Then
\[
\overline{{\mathfrak a} I^t \colon I^t} \subset \overline{{\mathfrak a}
I^t} \colon \overline{I}^t = (\overline{a}_s)\overline{I}^t \colon
\overline{I}^t \subset (\overline{a}_s) (\overline{I}^t \underset{Q}{\colon}
\overline{I}^t),
\]
where $Q$ is the total quotient ring of $\overline{R}$ and
$\overline{a}_s$ is a non zero divisor in $\overline{R}$ $($see
\cite[Proposition~3.1 and Proposition~3.3]{CEU}$)$.
By \cite[$(${\it c}$)$ in the proof of Theorem~4.1]{CEU}
we have that $(\overline{I}^{s-g})^{\vee\vee}$ is the
canonical module $\omega_{\overline{R}}$ of $\overline{R}$, where
$(\tratto)^{\vee} = {\rm Hom}_{\overline{R}}(\tratto, \omega_{\overline{R}})$.
Thus, for $1 \leq t \leq s-g$ it follows that
\[
\overline{I}^t \underset{Q}{\colon} \overline{I}^t \subset \overline{I}^{s-g}
\underset{Q}{\colon} \overline{I}^{s-g} \subset
(\overline{I}^{s-g})^{\vee\vee} \underset{Q}{\colon}
(\overline{I}^{s-g})^{\vee\vee} =
{\rm Hom}_{\overline{R}}(\omega_{\overline{R}},
\omega_{\overline{R}}) = \overline{R},
\]
as $\overline{R}$ satisfies property $S_2$ $($see \cite[$(${\it a}$)$ in the
proof of Theorem~4.1]{CEU}$)$. Thus $\overline{{\mathfrak a} I^t \colon I^t}
\subset (\overline{a}_s)$ and, after lifting back to the ring $R$, we get
\[
{\mathfrak a} I^t \colon I^t \subset (a_s) + K_{s-1} = {\mathfrak a}
+ K_{s-1}.
\]
Using Theorem~\ref{thm1}$(${\it b}$)$, for $1 \leq t \leq s-g$, we have that
\[
{\mathfrak a} I^t \colon I^t \subset ({\mathfrak a} + K_{s-1}) \cap I
= {\mathfrak a} + (K_{s-1} \cap I) = {\mathfrak a} + {\mathfrak a}_{s-1}
= {\mathfrak a},
\]
as desired. Again, we observe that the equality $K_{s-1} \cap I =
{\mathfrak a}_{s-1}$ follows from \cite[Proposition~3.4$(${\it c}$)$]{CEU}. \QED

\medskip

A case of interest is the following `global' formulation of
Theorem~\ref{thm2}.

\begin{Corollary}\label{nonlocal}
Let $R$ be a local Gorenstein ring. Let $I$ be an unmixed $R$-ideal of
grade $g$, satisfying $G_s$ for some $s \geq g$.
Suppose further that ${\rm Ext}^{g+j}_R(R/I^j, R)
= 0$ for $1 \leq j \leq s-g$. Let ${\mathfrak a} \colon I$ be an $s$-residual
intersection. Then for every $1 \leq t \leq s-g$
\[
{\mathfrak a} I^t \colon I^t = {\mathfrak a}.
\]
\end{Corollary}
\demo
The ideal $I$ is $(s-1)$-residually $S_2$ by \cite[Theorem~4.1]{CEU}.
The statement then follows from Theorem~\ref{thm2}. \QED

\medskip

\begin{Remark}
{\rm
\begin{itemize}
\item[$(${\it a}$)$]
The vanishing of the Ext modules in Theorem~\ref{thm2} and in
Corollary~\ref{nonlocal} are satisfied in particular if
${\rm depth}\, R/I^j \geq {\rm dim}\, R-g-j+1$ for
$1 \leq j \leq s-g$, which in turn holds if $I$ is a strongly Cohen--Macaulay
ideal $($assuming that $I$ satisfies $G_s)$.

\item[$(${\it b}$)$]
For applications to projective varieties, a convenient reformulation
of the vanishing of the Ext modules can be given in terms of vanishing of
certain local cohomology modules, via local duality $($see
\cite[Corollary~4.3]{CEU} for more details$)$.

\item[$(${\it c}$)$]
If we set $s=g+1$ then Corollary~\ref{nonlocal} recovers
\cite[Theorem~2.2]{H}.
\end{itemize} }
\end{Remark}

\subsection{The module case}

For the reader's sake, we start recalling some facts from \cite{RAM}. We will
always assume that $R$ is a local Gorenstein ring of dimension $d>0$,
infinite residue field and $E$ is a finitely generated $R$-module
with rank $e>0$. The {\em Rees algebra} ${\mathcal R}(E)$ of $E$ is the
symmetric algebra ${\mathcal S}(E)$ of $E$ modulo its $R$-torsion submodule.
If in addition $E$ is a submodule of a free $R$-module $G$, the above
definition coincide with the one, given by other authors, of the Rees
algebra of $E$ being the image of the natural map ${\mathcal S}(E)
\longrightarrow {\mathcal S}(G)$.

A submodule $U$ of $E$ is a {\em reduction} of $E$ or, equivalently,
$E$ is {\em integral} over $U$ if ${\mathcal R}(E)$ is integral over the
$R$-subalgebra generated by $U$.
Alternatively, the integrality condition is expressed by the equations
${\mathcal R}(E)_{r+1} = U \cdot {\mathcal R}(E)_r$ for $r \gg 0$.
The least integer $r \geq 0$ for which this equality holds is called the
{\em reduction number} of $E$ {\em with respect to} $U$ and denoted by
$r_U(E)$. The {\em reduction number} $r(E)$ of $E$ is defined to be the
minimum of $r_U(E)$, where $U$ ranges over all minimal reductions of $E$.
Since the residue field is assumed to be infinite, the minimal number of
generators of $U$ is given by the analytic spread $\ell=\ell(E)$ of $E$,
and it satisfies the inequalities: \ $e \leq \ell \leq d+e-1$.

In \cite{RAM} the study of the Rees algebra of a module $E$ is pursued via
the notion of a Bourbaki ideal of $E$. By a {\em Bourbaki ideal} of a
module $E$ we mean an ideal $I$ fitting into the short exact sequence
\[
0 \rightarrow F \longrightarrow E \longrightarrow I \rightarrow 0,
\]
with $F$ a free $R$-module. For technical reasons, it is actually better to
work with {\em generic} Bourbaki ideals, which are defined over a suitable
faithfully flat extension of $R$. Having this tool at our disposal, one can
then embark in the comparison of the Rees algebras of $E$ to the one of $I$.

Similar to what happens in the ideal case, also in many of the results
of \cite{RAM} a crucial role is played by conditions on the local number
of generators of a module. In the module case, $E$ is said to satisfy
condition $G_s$, for an integer $s \geq 1$, if $\mu (E_{\mathfrak p})\leq
{\rm dim}\, R_{\mathfrak p}+e-1$ whenever $1\leq {\rm dim}\,
R_{\mathfrak p}\leq s-1$.

\medskip

We are ready to state our main result in the case of modules. In the same
setting as Theorem~\ref{pd1}, it is worth pointing out that other
equivalent conditions to the Cohen--Macaulayness of ${\mathcal R}(E)$
are described in \cite[Theorem~4.7]{RAM}.

\begin{Theorem}\label{pd1}
Let $R$ be a local Gorenstein ring with infinite residue field and let
$E$ be a finitely generated $R$-module with ${\rm proj\, dim\,}E = 1$.
Write $e = {\rm rank\,} E$, $\ell = \ell(E)$ and assume that $E$ satisfies
$G_{\ell-e+1}$ and is torsionfree locally in codimension $1$.

If ${\mathcal R}(E)$ is Cohen--Macaulay,
then for any minimal reduction $U$ of $E$
\[
U \cdot {\mathcal R}(E)_t \, \underset{E}{\colon} \,
{\mathcal R}(E)_t = U
\]
for $0 \leq t \leq \ell - e - 1$.
\end{Theorem}
\demo
Write $U = \displaystyle\sum_{i=1}^{\ell} R a_i$. Let $z_{ij}$, with $1
\leq i
\leq \ell$ and $1 \leq j \leq e-1$, be variables and let
\[
\widetilde{R} = R(\{ z_{ij}\mbox{'s} \}) \qquad
\widetilde{E} = \widetilde{R} \otimes_R E \qquad
\widetilde{U} = \widetilde{R} \otimes_R U \qquad x_j = \sum_{i=1}^{\ell}
z_{ij}a_i \in \widetilde{U}.
\]
Consider $I \simeq \widetilde{E}/(x_1, \ldots, x_{e-1})$ and let
$J$ be the image of $\widetilde{U}$ in $I$. By \cite[Theorem~3.5]{RAM}
one has that
\begin{itemize}
\item
$I$ is a perfect ideal of grade $g=2$ $($hence {\em licci}$)$
satisfying $G_{\ell-e+1}$.

\item
Since $I$ is licci then ${\rm depth} \, R/I^j \geq {\rm dim} \, R -
g - j + 1$ for $1 \leq j \leq (\ell-e+1)-g+1$.

\item
Since $J$ is a reduction of $I$, one has that $J \colon I$ is an
$(\ell-e+1)$-residual intersection $($see \cite[Proposition~1.11]{U} or
\cite[Remark~2.7]{JU}$)$.

\item
${\mathcal R}(I)$ is Cohen--Macaulay, with ${\mathcal R}(I) \simeq
{\mathcal R}(\widetilde{E})/(x_1, \ldots, x_{e-1})$.
\end{itemize}
By Corollary~\ref{nonlocal} we conclude that
\[
JI^t \colon I^t = J
\]
for every $0 \leq t \leq (\ell-e+1)-g= \ell-e-1$.

Since ${\mathcal R}(I) \simeq {\mathcal R}(\widetilde{E})/(x_1,
\ldots, x_{e-1})$ we obtain that
\[
\widetilde{U} \subset
\widetilde{U} \cdot {\mathcal R}(\widetilde{E})_t \,
\underset{\widetilde{E}}{\colon} \, {\mathcal R}(\widetilde{E})_t
\subset \widetilde{U} + ((x_1, \ldots, x_{e-1}) \cdot
{\mathcal R}(\widetilde{E}))_1 \subset \widetilde{U},
\]
where the last inclusion holds since $((x_1, \ldots, x_{e-1}) \cdot
{\mathcal R}(\widetilde{E}))_1 \subset \widetilde{U}$.
Thus
\[
\widetilde{U} \cdot {\mathcal R}(\widetilde{E})_t
\, \underset{\widetilde{E}}{\colon} \, {\mathcal R}(\widetilde{E})_t
= \widetilde{U}.
\]
The final result now follows as $\widetilde{R}$ is a faithfully flat extension
of $R$. \QED

\medskip

\section{Applications to ideals}

The results in this Section recover
\cite[Corollary~2.3, Theorem~2.4, and Corollary~2.13]{H}.
Corollary~\ref{red} is a weaker version of
\cite[Corollary~5.5]{JU} and \cite[Corollary~2.4$(${\it c}$)$]{U2}.

\begin{Corollary}
Let $R$ be a local Gorenstein ring. Let $I$ be an $R$-ideal
of grade $g$, satisfying $G_s$ for some $s \geq g+1$ and unmixed
locally in codimension $s$.
Suppose further that $I$ is weakly $(s-1)$-residually $S_2$ and that
${\rm Ext}_R^{g+j}(R/I^j, R)_{\mathfrak p} = 0$
whenever ${\mathfrak p} \in V(I)$ with ${\rm dim}\, R_{\mathfrak p}\leq s$
and $1 \leq j \leq s - g$.
Let ${\mathfrak a} \colon I$ be an $s$-residual intersection and
suppose that ${\mathfrak a}$ and $I$ have the same radical $($e.g. if
$s$ is the analytic spread of $I$ and ${\mathfrak a}$ is a reduction of $I)$.
Then $I^n \subset {\mathfrak a}$ if and only if $I^{n+s-g} \cap
{\mathfrak a} \subset {\mathfrak a}I^{s-g}$.
\end{Corollary}
\demo
The proof is the one of \cite[Corollary~2.13]{H}.
If $I^n \subset {\mathfrak a}$ then $I^{n+s-g} = I^n I^{s-g} \subset
{\mathfrak a}I^{s-g}$. Conversely, suppose $I^{n+s-g} \cap {\mathfrak a}
\subset {\mathfrak a}I^{s-g}$. We know there exists $N$ such that $I^N
\subset {\mathfrak a}$. If $N \leq n$ then we are done. If $N-1 \geq n$, we
have $I^{N-1} I^{s-g} \subset I^{n+s-g} \cap {\mathfrak a} \subset
{\mathfrak a}I^{s-g}$ $($here we need $s-g\geq 1)$.
Thus $I^{N-1} \subset {\mathfrak a}I^{s-g} \colon
I^{s-g} = {\mathfrak a}$, by Theorem~\ref{thm2}.
Repeat the process until $I^n \subset {\mathfrak a}$.
\QED

\medskip

\begin{Corollary}\label{red}
Let $R$ be a local Gorenstein ring with infinite residue field.
Let $I$ be an $R$-ideal of grade $g$, satisfying $G_{\ell}$,
where $\ell$ is the analytic spread of $I$, and unmixed locally
in codimension $\ell$.
Suppose further that $I$ is weakly $(\ell-1)$-residually $S_2$ and that
${\rm Ext}_R^{g+j}(R/I^j, R)_{\mathfrak p} = 0$
whenever ${\mathfrak p} \in V(I)$ with ${\rm dim}\, R_{\mathfrak p}\leq \ell$
and $1 \leq j \leq \ell -g$. It follows that
\begin{itemize}
\item[$(${\it a}$)$]
the reduction number of $I$ is either zero or at least $\ell-g+1$;

\item[$(${\it b}$)$]
if $I^{r+1}=JI^r$, with $J$ a minimal reduction of $I$, then
$I^{r-\ell+g+1} \subset J$.
\end{itemize}
\end{Corollary}
\demo
$(${\it a}$)$
Assume the reduction number $r$ of $I$ to be different from zero and let
$J$ be a minimal reduction of $I$. Then $K = J\colon I$ is an
$\ell$-residual intersection $($see \cite[Proposition~1.11]{U} or
\cite[Remark~2.7]{JU}$)$. By Theorem~\ref{thm2}, we have that
\[
JI^t \colon I^t = J
\]
for every $1 \leq t \leq \ell-g$. Hence the statement.

$(${\it b}$)$
If $I^{r+1}=JI^r$ one has that $I^{r-\ell+g+1}I^{\ell-g} = JI^r \subset
JI^{\ell-g}$, as $r \geq \ell-g +1$ by $(${\it a}$)$. Thus
$I^{r-\ell+g+1} \subset JI^{\ell-g} \colon I^{\ell-g}=J$ by
Theorem~\ref{thm2}.
\QED

\medskip

In the next corollary, `$\,{}^{\tratto}$' denotes the integral closure
of an ideal.

\begin{Corollary}\label{intclos}
Let $R$ be a local regular ring with infinite residue field.
Let $I$ be a equidimensional $R$-ideal of grade $g$, satisfying $G_{\ell}$,
where $\ell$ is the analytic spread of $I$, and unmixed
locally in codimension $\ell$.
Suppose further that $I$ is weakly $(\ell-1)$-residually $S_2$ and that
${\rm Ext}_R^{g+j}(R/I^j,R)_{\mathfrak p} = 0$
whenever ${\mathfrak p} \in V(I)$ with ${\rm dim}\, R_{\mathfrak p}\leq \ell$
and $1 \leq j \leq \ell-g-1$.
If $I^{\ell-g}$ is unmixed locally in codimension $\ell$
then $\overline{I^g} \subset J$, for any reduction $J$ of $I$.
\end{Corollary}
\demo
We know that $I^{\ell-g} \overline{I^g} \subset \overline{I^{\ell}}
\subset JI^{\ell-g}$. The last inclusion follows from
\cite[Theorem~3.3]{AH2}, the fact that ${\rm bight}(I) = \max \{
{\rm ht}\, {\mathfrak p} \, : \, {\mathfrak p} \in {\rm min}(I) \} = g$ and
the assumption of $I^{\ell-g}$ being unmixed. Thus we conclude that
$\overline{I^g} \subset JI^{\ell-g} \colon I^{\ell-g} = J$,
by Theorem~\ref{thm2}. \QED


\medskip

\begin{thebibliography}{99}

\bibitem{AH2}{I.M. Aberbach and C. Huneke, An improved Brian\c{c}on--Skoda
theorem with applications to the Cohen--Macaulayness of Rees algebras,
Math. Ann. {\bf 297} (1993), 343--369.}

\bibitem{CEU}{M. Chardin, D. Eisenbud and B. Ulrich, Hilbert functions,
residual intersections, and residually $S_2$ ideals,
Compositio Math., to appear.}

\bibitem{CHV}{A. Corso, C. Huneke and W.V. Vasconcelos, On the
integral closure of ideals, Manuscripta Math. {\bf 95} (1998), 331--347.}

\bibitem{CPU}{A. Corso, C. Polini and B. Ulrich, Core of ideals and modules,
preprint 1999.}

\bibitem{H}{C. Huneke, A cancellation theorem for ideals,  J. Pure
\& Applied Algebra, to appear.}

\bibitem{JU}{M. Johnson and B. Ulrich, Artin--Nagata properties and
Cohen--Macaulay associated graded rings, Compositio Math. {\bf 103}
(1996), 7--29.}

\bibitem{LS}{J. Lipman and A. Sathaye,
Jacobian ideals and a theorem of Brian\c{c}on--Skoda, Michigan Math. J.
{\bf 28} (1981), 97--116.}

\bibitem{PU}{C. Polini and B. Ulrich, Linkage and reduction numbers,
Math. Ann. {\bf 310} (1998), 631--651.}

\bibitem{PU2}{C. Polini and B. Ulrich, Necessary and sufficient
conditions for the Cohen--Macaulayness of blowup algebras,
Compositio Math., to appear.}

\bibitem{RAM}{A. Simis, B. Ulrich and W.V. Vasconcelos, Rees
algebras of modules, preprint 1998.}

\bibitem{U}{B. Ulrich, Artin-Nagata properties and reductions of ideals,
in Commutative Algebra: Syzygies, Multiplicities, and Birational Algebra,
W. Heinzer, C. Huneke, J. Sally (eds.), Contemp. Math. {\bf 159}, Amer.
Math. Soc., Providence, 1994, 373--400.}

\bibitem{U2}{B. Ulrich, Ideals having the expected reduction number,
Amer. J. Math. {\bf 118} (1996), 17--38.}

\end{thebibliography}
\end{document}